\documentclass[12pt,leqno]{amsart}

\DeclareMathOperator{\rank}{rank}%
\DeclareMathOperator{\cd}{cd}%
\DeclareMathOperator{\cor}{cor}%
\DeclareMathOperator{\res}{res}%
\DeclareMathOperator{\ann}{ann}

\usepackage{hyperref}
\usepackage[arrow]{xypic}

\newcommand{\F}{\mathbb{F}}
\newcommand{\Fp}{\F_p}
\newcommand{\Gal}{\text{\rm Gal}}
\newcommand{\N}{\mathbb{N}}
\newcommand{\R}{\mathbb{R}}
\newcommand{\Z}{\mathbb{Z}}

\begin{document}

\title[Cohomological Dimension and Schreier's Formula]%
{Cohomological Dimension and Schreier's Formula in Galois Cohomology}

\begin{abstract}
    Let $p$ be a prime and $F$ a field containing a primitive $p$th
    root of unity.  Then for $n\in \N$, the cohomological dimension
    of the maximal pro-$p$-quotient $G$ of the absolute Galois group
    of $F$ is $\le n$ if and only if the corestriction maps $H^n(H,\Fp)
    \to H^n(G,\Fp)$ are surjective for all open subgroups $H$ of
    index $p$. Using this result we derive a surprising
    generalization to $\dim_{\Fp} H^n(H,\Fp)$ of Schreier's formula
    for $\dim_{\Fp} H^1(H,\Fp)$.
\end{abstract}

\author[Labute]{John Labute}
\address{Department of Mathematics and Statistics, \
McGill University, \linebreak Burnside Hall, 805 Sherbrooke Street
West, \ Montreal, Quebec \linebreak H3A 2K6 \ CANADA}
\email{labute@math.mcgill.ca}

\author[Lemire]{Nicole Lemire$^{\dag}$}
\address{Department of Mathematics, Middlesex College, \
University of Western Ontario, London, Ontario \ N6A 5B7 \ CANADA}
\thanks{$^\dag$Research supported in part by NSERC grant R3276A01.}
\email{nlemire@uwo.ca}

\author[Min\'{a}\v{c}]{J\'an Min\'a\v{c}$^{\ddag}$}
\thanks{$^\ddag$Research supported in part by NSERC grant R0370A01,
by the Mathematical Sciences Research Institute, Berkeley, and by a
2004/2005 Distinguished Research Professorship at the University of
Western Ontario.} \email{minac@uwo.ca}

\author[Swallow]{John Swallow}
\address{Department of Mathematics, Davidson College, Box 7046,
Davidson, North Carolina \ 28035-7046 \ USA}
\email{joswallow@davidson.edu}



\date{March 2, 2005}

\maketitle

\newtheorem{theorem}{Theorem}
\newtheorem*{proposition*}{Proposition}
\newtheorem{lemma}{Lemma}
\newtheorem*{corollary*}{Corollary}

\theoremstyle{definition}
\newtheorem*{remark*}{Remark}

\parskip=10pt plus 2pt minus 2pt

For a prime $p$, let $F(p)$ denote the maximal $p$-extension of a
field $F$.  One of the fundamental questions in the Galois theory
of $p$-extensions is to discover useful interpretations of the
cohomological dimension $\cd(G)$ of the Galois group
$G=\Gal(F(p)/F)$ in terms of the arithmetic of $p$-extensions of
$F$.  When $\cd(G)=1$, for instance, we know that $G$ is a free
pro-$p$-group \cite[\S 3.4]{S1}, and when $\cd(G)=2$ we have
important information on the $G$-module of relations in a minimal
presentation \cite[\S 7.3]{K}.

For a fixed $n > 2$, however, little is known about the structure
of $p$-extensions when $\cd(G)=n$.  Now when $n=1$ and $G$ is
finitely generated as a pro-$p$-group, we have Schreier's
well-known formula
\begin{equation}\label{eq:eq1}
    h_1(H) = 1 + [G:H](h_1(G) - 1)
\end{equation}
for each open subgroup $H$ of $G$, where
\begin{equation*}
        h_1(H) := \dim_{\Fp} H^1(H,\Fp).
\end{equation*}
(See, for instance, \cite[Example~6.3]{K}.)

Observe that from basic properties of $p$-groups it follows that
for each open subgroup $H$ of $G$ there exists a chain of
subgroups
\begin{equation*}
    G = G_0 \supset G_1 \supset \dots \supset
    G_k = H
\end{equation*}
such that $G_{i+1}$ is normal in $G_i$ and $[G_i:G_{i+1}]=p$ for
each $i=0,1,\dots,k-1$.  Since closed subgroups of free
pro-$p$-groups are free \cite[Corollary 3, \S I.4.2]{S1}, Schreier's
formula \eqref{eq:eq1} is equivalent to the seemingly weaker
statement that the formula holds for all open subgroups $H$ of $G$
of index $p$:
\begin{equation}\label{eq:eq2}
    h_1(H) = 1 + p(h_1(G) - 1).
\end{equation}

We deduce a remarkable generalization of Schreier's formula for
each $n\in \N$, as follows.  Let $F^\times$ denote the nonzero
elements of a field $F$, and for $c\in F^\times$, let $(c)\in
H^1(G,\Fp)$ denote the corresponding class.  For $\alpha\in
H^m(G,\Fp)$ abbreviate by $\ann_n \alpha$ the annihilator
\begin{equation*}
    \ann_n \alpha = \{ \beta \in H^n(G,\Fp) \ \ \vert \ \ \alpha
    \cup \beta = 0\}.
\end{equation*}
Finally, set $h_n(G)=\dim_{\Fp} H^n(G,\Fp)$.  Observe that the
hypothesis on the surjectivity of the corestriction on degree $n$
cohomology holds for all pro-$p$-groups of cohomological dimension
$n$ \cite[Proposition~3.3.8]{NSW}; conversely, in
section~\ref{se:when}, we show that this hypothesis for all
subgroups of index $p$ implies $\cd(G)\le n$.

\begin{theorem}\label{th:sf}
    Suppose that $\xi_p\in F$ and $h_n(G)<\infty$.  Let $H$ be an
    open subgroup of $G$ of index $p$, with fixed field
    $F(\root{p}\of{a})$, and suppose that the corestriction map
    $H^n(H,\Fp)\to H^n(G,\Fp)$ is surjective. Then
    \begin{equation*}
        h_n(H) = a_{n-1}(G,H) + p\big(h_n(G) - a_{n-1}(G,H)\big),
    \end{equation*}
    where $a_{n-1}(G,H)$ is the codimension of $\ann_{n-1} (a)$:
    \begin{equation*}
        a_{n-1}(G,H) := \dim_{\Fp} \big(\ H^{n-1}(G,\Fp)/ \ann_{n-1}
        (a)\ \big).
    \end{equation*}
\end{theorem}

The proof of Theorem~\ref{th:sf} brings additional insight into the
structure of Schreier's formula; in fact, it makes Schreier's
formula transparent for any $n\in \N$. In section~\ref{se:when}, we
derive several interpretations for the statement $\cd(G)=n$. First,
we prove in Theorem~\ref{th:s1t1} that if $F$ contains a primitive
$p$th root of unity $\xi_p$ then $\cd(G)\le n$ if and only if the
corestriction maps $\cor: H^n(H,\Fp) \to H^n(G,\Fp)$ are surjective
for all open subgroups $H$ of $G$ of index $p$.  As a corollary, we
show that the corresponding cohomology groups $H^{n+1}(H,\Fp)$ are
all free as $\Fp[G/H]$-modules if and only if $\cd(G)\le n$, under
the additional hypothesis that $F=F^2+F^2$ when $p=2$.  Finally, we
show in Theorem~\ref{th:s1t2} that if $G$ is finitely generated,
then $\cd(G)\le n$ if and only if a single corestriction map, from
the Frattini subgroup $\Phi(G)=G^p[G,G]$ of $G$, is surjective. In
section~\ref{se:sf} we prove Theorem~\ref{th:sf}.

For basic facts about Galois cohomology and maximal $p$-extensions
of fields, we refer to \cite{K} and \cite{S1}.  In particular, we
work in the category of pro-$p$-groups.

\section{When is $\cd(G)=n$?}\label{se:when}

As a consequence of recent results of Rost and Voevodsky on the
Bloch-Kato conjecture, we have the following interesting translation
of the statement $\cd(G) \le n$ for a given $n \in \N$.

\begin{theorem}\label{th:s1t1}
    Suppose that $\xi_p\in F$. Then for each $n \in \N$ we have
    $\cd(G)\le n$ if and only if
    \begin{equation*}
        \cor:H^n(H,\Fp)\to H^n(G,\Fp)
    \end{equation*}
    is surjective for every open subgroup $H$ of $G$ of index $p$.
\end{theorem}

\begin{proof}
    Suppose that $F$ satisfies the conditions of the theorem, and
    let $G_{F(p)}$ be the absolute Galois group of $F(p)$.

    Observe that since $F$ contains $\xi_p$, the maximal
    $p$-extension $F(p)$ is closed under taking $p$th roots and
    hence $H^1(G_{F(p)},\Fp)=\{0\}$.  By the Bloch-Kato
    conjecture, proved in \cite[Theorem~7.1]{V1}, the subring of
    the cohomology ring $H^\star(G_{F(p)}, \Fp)$ consisting of
    elements of positive degree is generated by cup-products of
    elements in $H^1(G_{F(p)},\Fp)$ .  Hence $H^n(G_{F(p)},\Fp) =
    \{0\}$ for $n\in \N$.  Then, considering the
    Lyndon-Hochschild-Serre spectral sequence associated to the
    exact sequence
    \begin{equation*}
        1 \to G_{F(p)} \to G_F \to G \to 1,
    \end{equation*}
    we have that
    \begin{equation}\label{eq:inf}
        \inf:H^\star(G,\Fp) \to H^\star(G_F,\Fp)
    \end{equation}
    is an isomorphism.

    Now suppose that $\cor : H^n(H,\Fp)\to H^n(G,\Fp)$ is surjective
    for all open subgroups $H$ of $G$ of index $p$.  Let $K$ be the
    fixed field of such a subgroup $H$. Then $K=F(\root{p}\of{a})$
    for some $a\in F^\times$. From Voevodsky's theorem
    \cite[Proposition~5.2]{V1}, modified in \cite[Theorem~5]{LMS1}
    and translated to $G$ from $G_F$ via the inflation maps
    \eqref{eq:inf} above, we obtain the following exact sequence:
    \begin{equation}\label{eq:es}
        H^n(H,\Fp) \xrightarrow{\cor} H^n(G,\Fp)
        \xrightarrow{\ -\cup(a)} H^{n+1}(G,\Fp)
        \xrightarrow{\res} H^{n+1}(H,\Fp).
    \end{equation}
    Therefore $\res:H^{n+1}(G,\Fp)\to H^{n+1}(H,\Fp)$ is injective
    for every open subgroup $H$ of $G$ of index $p$.

    Now consider an arbitrary element
    \begin{equation*}
        \alpha=(a_1)\cup\dots\cup(a_{n+1})\in H^{n+1}(G,\Fp),
    \end{equation*}
    where $a_i \in F^\times$ and $(a_i)$ is the element of
    $H^1(G,\Fp)$ associated to $a_i$, $i=1,2,\dots,n+1$. Suppose
    that $(a_1)\ne 0$, and set $K=F(\root{p}\of{a_1})$ and $H =
    \Gal(F(p)/K)$. We have $0=\res(\alpha)\in H^{n+1}(H,\Fp)$.
    Since $\res$ is injective, $\alpha=0$. Again by the Bloch-Kato
    conjecture \cite[Theorem~7.1]{V1}, we know that
    $H^{n+1}(G,\Fp)$ is generated by the elements $\alpha$ above.
    Hence $H^{n+1}(G,\Fp) = \{0\}$ and therefore $\cd(G)\le n$.
    (See \cite[page 49]{K}.)

    Conversely, if $\cd(G)\le n$ then by
    \cite[Proposition~3.3.8]{NSW} we conclude that $\cor:H^n(H,\Fp)
    \to H^n(G,\Fp)$ is surjective for open subgroups $H$ of $G$ of
    index $p$.
\end{proof}

Using conditions obtained in \cite{LMS2} for $H^n(H,\Fp)$ to be a
free $\Fp[G/H]$-module, we obtain the following corollary.  We
observe the convention that $\{0\}$ is a free $\Fp[G/H]$-module.

\begin{corollary*}
    Suppose that $\xi_p\in F$ and if $p=2$ suppose also that
    $F=F^2+F^2$. Then for each $n \in \N$, we have that
    $H^{n+1}(H,\Fp)$ is a free $\Fp[G/H]$-module for every open
    subgroup $H$ of $G$ of index $p$ if and only if $\cd(G)\le n$.
\end{corollary*}

Observe that the condition $F=F^2+F^2$ is satisfied in particular
when $F$ contains a primitive fourth root of unity $i$: for all
$c\in F^\times$, $c=((c+1)/2)^2 + ((c-1)i/2)^2$.

\begin{proof}
    Assume that $F$ is as above, $n\in \N$, and that $H^{n+1}
    (H,\Fp)$ is a free $\Fp[G/H]$-module for every open subgroup
    $H$ of $G$ of index $p$. If $p>2$, then it follows
    from \cite[Theorem 1]{LMS2} that the corestriction maps
    $\cor:H^n(H,\Fp)\to H^n(G,\Fp)$ are surjective for all such
    subgroups $H$.

    If $p=2$, then we consider open subgroups $H$ of index $2$
    with corresponding fixed fields $K=F(\sqrt{a})$.  From
    \cite[Theorem 1]{LMS2} we obtain that $\ann_n (a) = \ann_n
    \big( (a) \cup (-1)\big)$.  It follows from the hypothesis
    $F=F^2+F^2$ that $(c)\cup (-1) = 0 \in H^2(G,\F_2)$ for each
    $c\in F^\times$ and in particular for $c=a$.  Hence $\ann_n
    (a) = H^n(G,\F_2)$. But then from exact sequence \eqref{eq:es}
    above, we deduce that $\cor:H^n(H,\F_2)\to H^n(G,\F_2)$ is
    surjective.

    Since our analysis holds for all open subgroups
    $H$ of index $p$, by Theorem~\ref{th:s1t1} we conclude that
    $\cd(G)\le n$.

    Assume now that $\cd(G)\le n$. Then by Serre's theorem in
    \cite{S2} we find that $\cd(H)\le n$ for every open subgroup $H$
    of $G$.  Hence $H^{n+1}(H,\Fp)= \{0\}$ which, by our convention,
    is a free $\Fp[G/H]$-module, as required.
\end{proof}

\begin{remark*}
    When $p=2$ and $F\neq F^2 + F^2$, the statement of the
    corollary may fail. Consider the case $F=\R$. Then the only
    subgroup $H$ of index $2$ in $G=\Z/2\Z$ is $H=\{1\}$.  Then
    for all $n\in \N$, $H^{n+1}(H,\F_2) = \{0\}$ and is free as an
    $\F_2[G/H]$-module. However, $\cd(G)=\infty$.
\end{remark*}

Under the additional assumption that $G$ is finitely generated, we
show that the surjectivity of a single corestriction map is
equivalent to $\cd(G) \le n$.

\begin{theorem}\label{th:s1t2}
    Suppose that $\xi_p\in F$ and $G$ is finitely generated.  Then
    for each $n \in \N$ we have $\cd(G)\le n$ if and only if
    \begin{equation*}
        \cor:H^n(\Phi(G),\Fp)\to H^n(G,\Fp)
    \end{equation*}
    is surjective.
\end{theorem}

\begin{proof}
    Because $G$ is finitely generated, the index $[G:\Phi(G)]$
    is finite, and we may consider a suitable chain of open
    subgroups
    \begin{equation*}
        G = G_0 \supset G_1 \supset \dots \supset G_k =
        \Phi(G)
    \end{equation*}
    such that $[G_i:G_{i+1}]=p$ for each $i=0,1,\dots,k-1$.

    By Serre's theorem in \cite{S2}, $\cd(H)=\cd(G)$ for every
    open subgroup $H$ of $G$.  Hence if $\cd(G)\le n$ we may
    iteratively apply Theorem~\ref{th:s1t1} to the chain of open
    subgroups to conclude that
    \begin{equation*}
        \cor:H^n(\Phi(G),\Fp) \to H^n(G,\Fp)
    \end{equation*}
    is surjective.

    Assume now that $\cor:H^n(\Phi(G),\Fp)\to H^n(G,\Fp)$ is
    surjective. For each open subgroup $H$ of $G$ of index $p$
    we have a commutative diagram of corestriction maps
    \begin{equation*}
        \xymatrix{
        H^n(\Phi(G),\Fp) \ar[r] \ar[dr] &
        H^n(H,\Fp) \ar[d] \\ & H^n(G,\Fp) \\
        }
    \end{equation*}
    since $\Phi(G)\subset H$. We obtain that
    $\cor:H^n(H,\Fp)\to H^n(G,\Fp)$ is surjective, and by
    Theorem~\ref{th:s1t1} we deduce that $\cd(G) \le n$, as
    required.
\end{proof}

\section{Schreier's Formula for $H^n$}\label{se:sf}

We now prove Theorem~\ref{th:sf}.  Suppose that $H$ is an open
subgroup of $G$ of index $p$ and the corestriction map
$\cor:H^n(H,\Fp) \to H^n(G,\Fp)$ is surjective.  Let
$K=F(\root{p}\of{a})$ be the fixed field of $H$.

We claim that $\ann_{n-1} \big( (a)\cup (\xi_p) \big) =
H^{n-1}(G,\Fp)$.  Suppose that $\alpha\in H^{n-1}(G,\Fp)$.  By the
surjectivity hypothesis there exists $\beta\in H^n(H,\Fp)$ such that
$\cor \beta = (\xi_p) \cup \alpha$.  From Voevodsky's theorem
\cite[Proposition~5.2]{V1} modified in \cite[Theorem~5]{LMS1},
$(a)\cup (\cor\ \beta) = 0$ and hence $(a) \cup (\xi_p) \cup \alpha
= 0$.  Therefore the claim is established.

By \cite[Theorem~1]{LMS1}, we obtain the decomposition
\begin{equation*}
    H^n(H,\Fp) = X \oplus Y,
\end{equation*}
where $X$ is a trivial $\Fp[G/H]$-module and $Y$ is a free
$\Fp[G/H]$-module.  Moreover
\begin{align*}
    x &:= \dim_{\Fp} X = \dim_{\Fp} H^{n-1}(G,\Fp)/\ann_{n-1} (a)
     = a_{n-1}(G,H), \text{\ \ and}\\
    y &:= \rank\ Y = \dim_{\Fp} H^n(G,\Fp)/ (a) \cup
    H^{n-1}(G,\Fp).
\end{align*}
Therefore $h_n(H) = \dim_{\Fp} H^n(H,\Fp) = x + py$.

Now, considering the exact sequence
\begin{equation*}
    0 \to \frac{H^{n-1}(G,\Fp)}{\ann_{n-1}(a)} \xrightarrow{\ -\cup
    (a)} H^n(G,\Fp) \to \frac{H^n(G,\Fp)} {(a)\cup H^{n-1}(G,\Fp)}
    \to 0,
\end{equation*}
we see that $\dim_{\Fp} H^n(H,\Fp)$ is equal to the sum of the
dimension $x$ of the kernel and $p$ times the dimension $y$ of the
cokernel, and the theorem follows.

Observe that we have established a more general formula than the
formula displayed in Theorem~\ref{th:sf}, since we have not
assumed that $h_n(G)$ is finite.

When $n=1$, $\ann_{n-1} (a) = \{0\}$ so that $a_{n-1}(G,H)=1$.
Therefore when $G$ is finitely generated we recover Schreier's
formula \eqref{eq:eq2}:
\begin{equation*}
    h_1(H)=1 + p(h_1(G)-1).
\end{equation*}

\end{document}